\newif\ifprint 
\printfalse 

\newif\ifdraft
\newif\ifonly
\newif\ifsuppress
\newif\ifhyper
\newif\ifproceedings
\newif\ifarxiv 

\arxivfalse
\suppresstrue
\hypertrue

\newif\ifshowprivate 
\newif\ifbells 

\showprivatefalse
\bellstrue

\ifsuppress
	\gdef\wlog#1{}
	\global\let\GenericInfo\relax
	\DeclareRobustCommand{\GenericInfo}[2]{}
\fi

\def\commentx#1{}
\commentx{
\documentclass{mrlart4}

  \setcounter{page}{10001}
  \overfullrule=5pt
}

\documentclass[10pt]{amsart}

\draftfalse
\ifdraft
	\calclayout
\fi

\usepackage{amssymb,latexsym,amsmath,amsthm}
	
\usepackage{xspace}
\usepackage{color}
\usepackage{diagrams} 
\usepackage[pdftex]{graphicx}

\ifdraft
	\def\pdfsyncstop{}
	\def\pdfsyncstart{}
\else
\def\pdfsyncstop{}
\def\pdfsyncstart{}
\fi

\ifhyper
\ifprint
	\definecolor{linkred}{rgb}{0,0,0} 
	\definecolor{linkblue}{rgb}{0,0,0} 
\else
	\definecolor{linkred}{rgb}{0.7,0.2,0.2}
	\definecolor{linkblue}{rgb}{0,0.2,0.6}
\fi
\usepackage[
	hyperindex,
	pagebackref,
	pdftex,
	pdfusetitle,
	breaklinks=true,
	bookmarks=false,
	plainpages=false,
	colorlinks,
	linkcolor=linkblue,
	citecolor=linkblue,
	urlcolor=linkred,
	pdftoolbar=false,
	pdfpagelabels,
	pdfstartview=Fit
]
{hyperref}
\usepackage{memhfixc}
\usepackage[backrefs,msc-links,nobysame]{amsrefs}
\usepackage[low-sup]{subdepth}

\newtheorem{theorem}{Theorem}
\newtheorem{lemma}[theorem]{Lemma}
\newtheorem{proposition}[theorem]{Proposition}
\newtheorem{corollary}[theorem]{Corollary}

\theoremstyle{definition}

\theoremstyle{remark}

\def\mcl{\mathcal}

\def\ten{\otimes}

\def\b{\beta}
\def\d{\delta}

\def\o{\omega}

\def\SS{\mcl{O}}
\def\deg{\textup{deg} \, }

\def\Hg2{Hilb_{g,2}}
\def\Cg2{Chow_{g,2}}

\def\om2{\omega^{\ten 2}}

\def\Hilb{\text{Hilb}}
\def\Chow{\text{Chow}}
\def\chowq{\Chow_{g,4}/\!\!/SL_{7g-7}}
\def\hilbq{\Hilb_{g,4}/\!\!/SL_{7g-7}}

\def\k{k}
\def\resp#1{{\upshape{[}}resp. #1{\upshape{]}}}

\def\setdrawbox#1#2#3{
\pdfximage{#2} 
\setbox0=\hbox{\pdfrefximage\pdflastximage} 
\drawx=#1\wd0
\ifdim\drawx>\hsize\drawx=\hsize\fi
\pdfximage width \drawx {#2}
\setbox\drawing=\vbox{\offinterlineskip\pdfrefximage\pdflastximage\kern 0pt}
\drawx=\wd\drawing
\drawy=\ht\drawing
\ngap=0pt \sgap=0pt \wgap=0pt \egap=0pt
\setbox0=\vbox{\offinterlineskip \box\drawing \ifgridlines\drawgrid\drawx\drawy\fi #3}%
\setbox\drawing=\vbox{\kern\ngap\hbox{\kern\wgap\box0\kern\egap}\kern\sgap}
}
\def\obsoletedrawbox#1#2#3{\vbox{
  \setbox\drawing=\vbox{\offinterlineskip\epsfbox{#2.eps}\kern 0pt}
  \drawbp=\epsfurx
  \advance\drawbp by-\epsfllx\relax
  \multiply\drawbp by #1
  \divide\drawbp by 100
  \drawx=\drawbp truebp
  \ifdim\drawx>\hsize\drawx=\hsize\fi
  \epsfxsize=\drawx
  \setbox\drawing=\vbox{\offinterlineskip\epsfbox{#2.eps}\kern 0pt}
  \drawx=\wd\drawing
  \drawy=\ht\drawing
  \ngap=0pt \sgap=0pt \wgap=0pt \egap=0pt
  \setbox0=\vbox{\offinterlineskip
    \box\drawing \ifgridlines\drawgrid\drawx\drawy\fi #3}
  \kern\ngap\hbox{\kern\wgap\box0\kern\egap}\kern\sgap}}

\def\drawsmash#1{\relax
\ifmmode%
\typeout{Warning:Macro drawsmash used in math mode}%
\fi%
\setbox0=\hbox{#1}\ht0=0pt\dp0=0pt\box0}

\newcount\figcount
\newbox\drawing
\newcount\drawbp
\newdimen\drawx
\newdimen\drawy
\newdimen\ngap
\newdimen\sgap
\newdimen\wgap
\newdimen\egap

\newif\ifgridlines
\newbox\figtbox
\newbox\figgbox
\newdimen\figtx
\newdimen\figty

\newdimen\bwd
\def\hhline#1{\vbox{\drawsmash{\hbox to #1{\leaders\hrule height \bwd\hfil}}}}
\def\vvline#1{\hbox to 0pt{%
  \hss\vbox to #1{\leaders\vrule width \bwd\vfil}\hss}}

\def\clap#1{\hbox to 0pt{\hss#1\hss}}
\def\vclap#1{\vbox to 0pt{\offinterlineskip\vss#1\vss}}

\def\hstutter#1#2{\hbox{%
  \setbox0=\hbox{#1}%
  \hbox to #2\wd0{\leaders\box0\hfil}}}

\def\vstutter#1#2{\vbox{
  \setbox0=\vbox{\offinterlineskip #1}
  \dp0=0pt
  \vbox to #2\ht0{\leaders\box0\vfil}}}

\def\crosshairs#1#2{
  \dimen1=.002\drawx
  \dimen2=.002\drawy
  \ifdim\dimen1<\dimen2\dimen3\dimen1\else\dimen3\dimen2\fi
  \setbox1=\vclap{\vvline{2\dimen3}}
  \setbox2=\clap{\hhline{2\dimen3}}
  \setbox3=\hstutter{\kern\dimen1\box1}{4}
  \setbox4=\vstutter{\kern\dimen2\box2}{4}
  \setbox1=\vclap{\vvline{4\dimen3}}
  \setbox2=\clap{\hhline{4\dimen3}}
  \setbox5=\clap{\copy1\hstutter{\box3\kern\dimen1\box1}{6}}
  \setbox6=\vclap{\copy2\vstutter{\box4\kern\dimen2\box2}{6}}
  \setbox1=\vbox{\offinterlineskip\box5\box6}
  \drawsmash{\vbox to #2{\hbox to #1{\hss\box1}\vss}}}

\def\boxgrid#1{\rlap{\vbox{\offinterlineskip
  \setbox0=\hhline{\wd#1}
  \setbox1=\vvline{\ht#1}
  \drawsmash{\vbox to \ht#1{\offinterlineskip\copy0\vfil\box0}}
  \drawsmash{\vbox{\hbox to \wd#1{\copy1\hfil\box1}}}}}}

\def\drawgrid#1#2{\vbox{\offinterlineskip
  \dimen0=\drawx
  \dimen1=\drawy
  \divide\dimen0 by 10
  \divide\dimen1 by 10
  \setbox0=\hhline\drawx
  \setbox1=\vvline\drawy
  \drawsmash{\vbox{\offinterlineskip
    \copy0\vstutter{\kern\dimen1\box0}{10}}}
  \drawsmash{\hbox{\copy1\hstutter{\kern\dimen0\box1}{10}}}}}

\def\figtext#1#2#3#4#5{
  \setbox\figtbox=\hbox{#5}
  \dp\figtbox=0pt
  \figtx=-#3\wd\figtbox \figty=-#4\ht\figtbox
  \advance\figtx by #1\drawx \advance\figty by #2\drawy
  \dimen0=\figtx \advance\dimen0 by\wd\figtbox \advance\dimen0 by-\drawx
  \ifdim\dimen0>\egap\global\egap=\dimen0\fi
  \dimen0=\figty \advance\dimen0 by\ht\figtbox \advance\dimen0 by-\drawy
  \ifdim\dimen0>\ngap\global\ngap=\dimen0\fi
  \dimen0=-\figtx
  \ifdim\dimen0>\wgap\global\wgap=\dimen0\fi
  \dimen0=-\figty
  \ifdim\dimen0>\sgap\global\sgap=\dimen0\fi
  \drawsmash{\rlap{\vbox{\offinterlineskip
    \hbox{\hbox to \figtx{}\ifgridlines\boxgrid\figtbox\fi\box\figtbox}
    \vbox to \figty{}
    \ifgridlines\crosshairs{#1\drawx}{#2\drawy}\fi
    \kern 0pt}}}}

\def\nwtext#1#2#3{\figtext{#1}{#2}01{#3}}
\def\netext#1#2#3{\figtext{#1}{#2}11{#3}}
\def\swtext#1#2#3{\figtext{#1}{#2}00{#3}}

\def\wtext#1#2#3{\figtext{#1}{#2}0{.5}{#3}}
\def\etext#1#2#3{\figtext{#1}{#2}1{.5}{#3}}
\def\ntext#1#2#3{\figtext{#1}{#2}{.5}1{#3}}
\def\stext#1#2#3{\figtext{#1}{#2}{.5}0{#3}}

\def\hpad#1#2#3{\hbox{\kern #1\hbox{#3}\kern #2}}
\def\vpad#1#2#3{\setbox0=\hbox{#3}\dp0=0pt\vbox{\kern #1\box0\kern #2}}

\def\stack#1#2#3{\vbox{\offinterlineskip
  \setbox2=\hbox{#2}
  \setbox3=\hbox{#3}
  \dimen0=\ifdim\wd2>\wd3\wd2\else\wd3\fi
  \hbox to \dimen0{\hss\box2\hss}
  \kern #1
  \hbox to \dimen0{\hss\box3\hss}}}

\def\hexp#1{%
  \setbox0=\hbox{${}^{#1}$}%
  \hbox to .5\wd0{\box0\hss}}

\long\def\commentout#1{}

\renewcommand{\eprint}[1]{#1}

\makeatletter

\def\bib@div@mark#1{%
    \@mkboth{{#1}}{{#1}}%
	}
\def\print@backrefs#1{%
    \space\SentenceSpace$\leftarrow$\csname br@#1\endcsname
}
\catcode`\'=11 
\renewcommand{\PrintAuthors}[1]{%
    \ifx\previous@primary\current@primary
        \sameauthors\@empty
    \else
        \def\current@bibfield{\bib'author}%
		        \PrintNames{}{}{\scshape #1}%
   \fi
}
\catcode`\'=12
\ifarxiv
\else
\def\MRhref#1#2{%
    \begingroup
        \parse@MR#1 ()\@empty\@nil%
        \href{\MR@url}{\texttt{\@tempd\vphantom{()}}}%
        \ifx\@tempe\@empty
        \else
            \ \href{\MR@url}{\texttt{(\@tempe)}}%
        \fi
    \endgroup
}%
\def\MR#1{%
    \relax\ifhmode\unskip\spacefactor3000 \space\fi
    \begingroup
        \strip@MRprefix#1\@nil
        \edef\@tempa{\@nx\MRhref{MR\@tempa}{\@tempa}}%
    \@xp\endgroup
    \@tempa
}
\fi

\makeatother

\newcommand\capdraw[5]{
\begin{figure}[tb]
  \setdrawbox{#3}{#4}{#5}
  \centerline{\ifgridlines\boxgrid\drawing\fi\box\drawing}
  \caption{#2}
  \label{#1}
\end{figure}%
}

\newarrow{Injectto}{>}{-}{-}{-}{>}
\newarrow{Identifyto}{=}{=}{=}{=}{=}
\newarrow{Rationalto}{}{dash}{}{dash}{>}
\newarrow{IMapsto}{|}{-}{-}{-}{>}

\newdiagramgrid{ilgrid}{0.385, 0.385}{}
\def\ilto{\pdfsyncstop\mathop{\hbox{\kern -5pt
\begin{diagram}[midshaft,inline,grid=ilgrid]
&\rTo&\end{diagram}\kern -5pt}}\pdfsyncstart}
\def\ilinto{\pdfsyncstop\mathop{\hbox{\kern -5pt
\begin{diagram}[midshaft,inline,grid=ilgrid]
&\rInjectto&\end{diagram}\kern -5pt}}\pdfsyncstart}
\def\ilmapsto{\pdfsyncstop\mathop{\hbox{\kern -5pt
\begin{diagram}[midshaft,inline,grid=ilgrid]
&\rIMapsto&\end{diagram}\kern -5pt}}\pdfsyncstart}
\newdiagramgrid{ildashgrid}{0.55, 0.55}{}
\def\ildashto{\pdfsyncstop\mathop{\hbox{\kern -5pt
\begin{diagram}[midshaft,inline,grid=ildashgrid]
&\rDashto&\end{diagram}\kern -5pt}}\pdfsyncstart}

\renewcommand{\to}{\pdfsyncstop\ilto\pdfsyncstart}
\renewcommand{\mapsto}{\pdfsyncstop\ilmapsto\pdfsyncstart}

\newdiagramgrid{esgrid}{0.585, 0.585}{}
\def\esarrowplain{\pdfsyncstop\mathop{\hbox{\kern -3pt
\begin{diagram}[midshaft,inline,grid=esgrid]
&\rTo&\end{diagram}\kern -3pt}}\pdfsyncstart}
\def\esarrowmap#1{\pdfsyncstop\mathop{\hbox{\kern -3pt
\begin{diagram}[midshaft,inline,grid=esgrid]
&\rTo&\end{diagram}\kern -3pt}}^{\kern -3pt #1}\pdfsyncstart}
\def\esarrowinject{\pdfsyncstop\mathop{\hbox{\kern -3pt
\begin{diagram}[midshaft,inline,grid=esgrid]
&\rInjectto&\end{diagram}\kern -3pt}}\pdfsyncstart}
\def\esarrowinjectmap#1{\pdfsyncstop\mathop{\hbox{\kern -3pt
\begin{diagram}[midshaft,inline,grid=esgrid]
&\rInjectto&\end{diagram}\kern -3pt}}^{#1}\pdfsyncstart}
\newdiagramgrid{esgridshort}{0.47, 0.47}{}
\def\shortesarrowplain{\pdfsyncstop\mathop{\hbox{\kern -3pt
\begin{diagram}[midshaft,inline,grid=esgridshort]
&\rTo&\end{diagram}\kern -3pt}}\pdfsyncstart}
\def\shortesarrowmap#1{\pdfsyncstop\mathop{\hbox{\kern -3pt
\begin{diagram}[midshaft,inline,grid=esgridshort]
&\rTo&\end{diagram}\kern -3pt}}^{\kern -3pt #1}\pdfsyncstart}
\def\shortesarrowinject{\pdfsyncstop\mathop{\hbox{\kern -3pt
\begin{diagram}[midshaft,inline,grid=esgridshort]
&\rInjectto&\end{diagram}\kern -3pt}}\pdfsyncstart}
\def\shortesarrowinjectmap#1{\pdfsyncstop\mathop{\hbox{\kern -3pt
\begin{diagram}[midshaft,inline,grid=esgridshort]
&\rInjectto&\end{diagram}\kern -3pt}}^{#1}\pdfsyncstart}

\newcommand{\setdiagram}[3]{
\vskip-12pt
\ifx#1\endA\endB
\setbox0=\hbox{}
\else
\refstepcounter{equation}
\setbox0=\hbox{\thmandnumberfont(\theequation)}
\label{#1}
\fi
\pdfsyncstop
\ifx#2\endA\endB
\begin{diagram}[eqno=\box0,midshaft]
#3
\end{diagram}
\else
\begin{diagram}[eqno=\box0,midshaft,#2]
#3
\end{diagram}
\fi
\pdfsyncstart
}

\newcommand{\inlinechoose}[2]%
{\smash{\binom{#1}{#2}}}

\newcommand{\inlinefrac}[2]%
{\smash{\bigl({\frac{#1}{#2}} \bigr)}}

\newcommand{\inlinefracnoparen}[2]%
{\smash{ {\frac{#1}{#2}} }}

\newcommand{\PP}{\ensuremath{\mathbb{P}}}

\newcommand{\spec}{\mathrm{Spec}} 

\newcommand{\ssl}{\mathrm{SL}}

\newcommand{\diag}{\mathrm{diag}}

\newcommand{\shavedast}{\ast \kern -1.25pt}

\renewcommand{\epsilon}{\varepsilon}

\newcommand{\thst}[2]{\ensuremath{{#1}^{\mathrm{#2}}}}

 \newcommand{\restrictedto}[2]{%
 \ensuremath{{#1\hskip 1pt{\vrule height 7.2pt depth 3.6pt}\hskip0.75pt}%
 \raisebox{-.8pt}[0pt][0pt]{\ensuremath{\null_{#2}}}}}

\newcommand{\h}{\ensuremath{H}}
\newcommand{\m}{\ensuremath{M}}

\renewcommand{\o}{\ensuremath{\mathcal{O}}} 
\renewcommand{\i}{\ensuremath{\mathcal{I}}} 

\renewcommand{\phi}{\varphi}

\newcommand{\mbar}{\ensuremath{\overline{\m}}}
\newcommand{\mibar}[1]{\ensuremath{\mbar_{#1}}}

\newcommand{\mgbar}{\ensuremath{\mibar{g}}}

\renewcommand{\k}{\ensuremath{\mathcal{K}}}

\renewcommand{\hbar}{\ensuremath{\overline{H}}}

\newcommand{\semi}[1]{\ensuremath{#1^{\text{ss}}}}

\newcommand{\ord}{\mathrm{ord}}

\newcommand{\invv}[2]{{#1}^{\kern #2{\null^{-1}}}}

\newcommand{\GIT}{GIT\xspace}

\newcommand{\init}{\text{in}}

\newcommand{\lbpow}[2]{#1^{\otimes #2}}
\newcommand{\lbm}[1]{#1^{\otimes m}}

\newcommand{\HO}[2]{H^0\bigl(#1,#2\bigr)}
\newcommand{\hO}[2]{h^0\bigl(#1,#2\bigr)}

\newcommand{\mpsg}{\overline{M}_g^{\text{ps}}}

\newcommand{\flushsubsection}[1]{\subsection{\kern-\normalparindent #1}}
\newcommand{\flushsubsectionstar}[1]{\subsection*{\kern-\normalparindent #1}}

\normalfont\selectfont
\begin{document}

\title{Stability of Tails and $4$-Canonical Models}

\author{Donghoon Hyeon}
\address{Department of Mathematics\\ Marshall University\\ Huntington, WV 25755}
\email{hyeond@marshall.edu}

\author{Ian Morrison}
\address{Department of Mathematics\\ Fordham University\\ Bronx, NY 10458}
\curraddr{Mathematical Sciences Research Institute\\ Berkeley CA 94720}
\email{morrison@fordham.edu}


\begin{abstract}
	We show that the GIT quotients of suitable loci in the Hilbert and Chow schemes of $4$-canonically embedded curves of genus $g\ge 3$ are the moduli space $\overline{M}_g^{\text{ps}}$ of pseudostable curves constructed by Schubert in~\cite{Schubert} using Chow varieties and $3$-canonical models. The only new ingredient needed in the Hilbert scheme variant is a more careful analysis of the stability with respect to a certain $1$-ps $\lambda$ of the $m^{\text{th}}$ Hilbert points of curves $X$ with elliptic tails. We compute the exact weight with which $\lambda$ acts, and not just the leading term in $m$ of this weight. A similar analysis of stability of curves with rational cuspidal tails allows us to determine the stable and semistable $4$-canonical Chow loci. Although here the geometry of the quotient is more complicated because there are strictly semistable orbits, we are able to again identify it as $\overline{M}_g^{\text{ps}}$. Our computations yield, as byproducts, examples of both $m$-Hilbert unstable and $m$-Hilbert stable $X$ that are Chow strictly semistable.
\end{abstract} 

\maketitle

\subsection*{Introduction} Over 30 years ago, Gieseker~\cites{GiesekerTata, GiesekerCIME} gave a  GIT construction of the moduli space $\mgbar$ of stable curves $X$ of genus $g$ as a quotient of a suitable locus in the Hilbert scheme of $\nu$-canonically embedded curves and Mumford~\cite{MumfordEnseignement} gave a variant using Chow points. A condensed exposition is given in~\cite{Moduli}*{Chapter 4}. These constructions require that $\nu \ge 5$. Since $\lbpow{\omega_X}{\nu}$ is very ample when $\nu \ge 3$, a natural question is to describe the quotient when $\nu$ equals $3$ or $4$.

The only point in Gieseker's construction at which the hypothesis $\nu \ge 5$ cannot be weakened is in applying the calculation that ordinary cusps (singularities locally analytically isomorphic to $y^2=x^3$) destabilize Hilbert and Chow points of curves of degree $d$ in $\PP^{n-1}$ when $\frac{d}{n}< \frac{8}{7}$: for $\nu$-canonical curves $\frac{d}{n}=\frac{2v}{2v-1}$. For a proof with this constant, which is sharp as we will see in Lemma~\ref{cusplambdainverse}, see \cite{MumfordEnseignement}*{Proposition~3.1}; most references---\cite{GiesekerTata}*{Proposition~1.0.5}, \cite{GiesekerCIME}*{Lemma~5.7} and \cite{Moduli}*{Step 4 on p. 251}---give the weaker bound $\frac{9}{8}$ which is easier to verify and suffices for Gieseker's application. In view of this result, $\nu \ge 5$  is needed implicitly in the proof, by semistable replacement, that all reducible stable curves are Hilbert stable.  

Recall that an elliptic tail is a connected genus $1$ subcurve of a stable curve meeting the rest of the curve at a single point and a varying elliptic tail is a curve $\PP^1 \subset \mgbar$ obtained by gluing an elliptic tail of varying $j$-invariant to a fixed point on a fixed curve of genus $g-1$. Recall also that a divisor on $\mgbar$ is said to have slope $s$ if it is a positive multiple of $s\lambda - \delta$. On the one hand (cf.\cite{Moduli}*{\S 3.F}), divisors of slope $11$ are numerically trivial on varying elliptic tails. On the other, Mumford~\cite{MumfordEnseignement}*{Section~5} (see also~\cites{HassettHyeonFlip, MorrisonGIT}) showed that, exactly when $\nu \le 4$, the polarizations with which \GIT naturally equips $\nu$-canonical quotients have slope $s \le 11$.

These observations suggest that ordinary cusps appear and elliptic tails disappear in the $\nu$-canonical quotients for $\nu <5$. This prediction was verified by Schubert~\cites{Schubert, MorrisonGIT} who, using $3$-canonical Chow points, produced, for $g \ge 3$, a quotient that is a coarse moduli space $\mpsg$ for pseudostable curves: these are reduced, connected, complete curves with finite automorphism group, whose only singularities are nodes and ordinary cusps, and which have no elliptic tails. If $g=2$, the subcurve $C$ in \cite{Schubert}*{Lemma 4.1} is not unique and the moduli problem is not separated but this quotient was later constructed by Hyeon and Lee \cite{HyeonLeeGenusThree} who showed that all pseudostable cuspidal curves are identified in it. Recently Hassett and Hyeon~\cite{HassettHyeonFlip} have analyzed the case $\nu=2$ producing moduli spaces in which cusps and tacnodes are allowed but elliptic tails and certain elliptic chains are excluded. But the case $\nu=4$, whose GIT you'd expect to be easier, has never been treated. The aim of this short note is to fill in this gap and confirm the expectation that the quotient is again $\mpsg$.  We work throughout over an algebraically closed field $k$ of arbitrary characteristic.

A careful reading of Schubert's paper shows that, with one exception, all the instability and stability claims made in Sections~1--3 for $\nu=3$ remain valid when $\nu=4$ with only minor numerical modifications to reflect the change in $\nu$. As we will see in Corollary~\ref{fourcanonicaltailsunstable}, his proof that curves with elliptic tails have Chow unstable $3$-canonical models only implies that their $4$-canonical models are not Chow stable.

Schubert's argument can be rescued in two ways. The first is to work with Hilbert points. Indeed, his construction applies, mutatis mutandis, using Hilbert points of $3$-canonical curves; there are no Chow strictly semistable points and Chow stability or instability implies the corresponding Hilbert property (cf. \cite{HassettHyeonFlip}*{Proposition 3.13}). To make the argument apply to Hilbert points of $4$-canonical curves requires only a more precise analysis of the weights with which an explicit one-parameter subgroup (henceforth $1$-ps) $\rho$ acts on the $\thst{m}{th}$ Hilbert points of certain projective embeddings of stable curves with elliptic tails. It would suffice to bound below the quadratic and linear terms in $m$ of these weights, but, we are able to compute them exactly in Lemma~\ref{lambdaweightlemma}. This calculation leads easily to the conclusion in Corollary~\ref{fourcanonicaltailsunstable} that $\thst{m}{th}$ Hilbert points of $4$-canonical models of stable curves with elliptic tails are $\rho$-unstable. 

Substituting Corollary~\ref{fourcanonicaltailsunstable} for his Lemma~3.1, the remainder of Schubert's construction, with the Chow scheme for $3$-canonical curves replaced by the  $\thst{m}{th}$ Hilbert scheme of $4$-canonical curves for a sufficiently large $m$, goes through with only the sort of numerical changes mentioned.

This leaves open the question of the $4$-canonical Chow quotient. Here the quotient again turns out to be $\mpsg$ but the geometry is a bit more complicated. The closures of the Chow orbits both of  cuspidal curves and of curves with elliptic tails contain the orbits of curves with (rational) cuspidal tails and all three types are Chow strictly semistable. To see that none are stable, we first compute the weight with which the same $\lambda$ acts on $4$-canonical curves with cuspidal tails by an even more direct argument in Lemma~\ref{L:Cstar}. This allows us to compute basins of attraction for this $\lambda$ and its inverse in Lemma~\ref{L:basin}. A semistable replacement argument (Proposition~\ref{P:chowfour}) then shows that these orbits are not unstable. This lets us, in Theorem~\ref{T:main}, use the cycle map to see that the Chow and Hilbert quotients are isomorphic and hence identify the former with $\mpsg$.

Corollary~\ref{fourcanonicaltailsunstable} also shows that the $4$-canonical Chow points of curves $X$ with elliptic tails are strictly semistable and Hilbert unstable with respect to $\rho$ and Lemma~\ref{L:Cstar} gives the same conclusions for $Y$ with cuspidal tails. Our $4$-canonical Chow quotient construction shows that such curves are Chow strictly semistable. Likewise, as a complement to Lemma~\ref{L:basin}, we compute in Lemma~\ref{cusplambdainverse} the weights with which  $\rho^{-1}$ acts on a $4$-canonical cuspidal curve $Z$ and show that such curves are Chow strictly semistable and $m$-Hilbert stable with respect to $\rho^{-1}$. The same conclusions follow unconditionally (i.e. for any $1$-ps) for such $Z$ from the constructions of the $4$-canonical Chow and Hilbert quotients.

Such $X$ and $Y$ are the easiest examples known to us of projective embeddings with strictly semistable Chow point and unstable $m$-Hilbert points and such $Z$ the only ones with strictly semistable Chow point and stable $m$-Hilbert points. Both sets of examples are, in a sense, optimal. Corollary~\ref{fourcanonicaltailsunstable} and Lemma~\ref{L:Cstar} show that, with respect to a critical $1$-ps, their Hilbert-Mumford indices are, respectively, $-(m-1)$ and $(m-1)$ and \cite{HassettHyeonFlip}*{Proposition 3.17} implies that these indices must be divisible by $(m-1)$.\footnote{~This divisibility is used in proving Lemma~\ref{L:Cstar}, but it gives a check on the calculations of Corollary~\ref{fourcanonicaltailsunstable} and Lemma~\ref{cusplambdainverse}.} Other examples of the former type arise in the analysis of bicanonically embedded rosaries in~\cite{HassettHyeonFlip}*{Propositions~10.1 and~10.7}. Any Veronese embedding of a projective space is both Chow and Hilbert strictly semistable with respect to any $\rho$ lying in its stabilizer, and also unconditionally.

\subsection*{Acknowledgements} The impetus for this paper came from Dave Swinarski who noted that Schubert~\cite{Schubert} treats only $\nu=3$ and asked what happens for $\nu=4$. The first author was supported by Korea Science and Engineering Foundation (KOSEF) grant No. R01-2007-000-10948-0. The second author was supported by a Fordham University Faculty Fellowship and the paper was written while he was visiting the University of Sydney; he thanks both institutions and Gus Lehrer, Amnon Neeman and Paul Norbury for their hospitality during his time in Australia.  

\subsection*{The $1$-ps $\rho$ and stability with respect to it}
If $X$ is a Deligne-Mumford stable curve that can be written $X = C \cup E$ where $C$ and $E$ are subcurves of genera $(g-1)$ and $1$ respectively and $C\cap E$ is a single node $p$, then we say that $X$ has an elliptic tail. We do not assume $C$ smooth or irreducible. Let $L$ be a very ample, non-special line bundle on $X$ of degree $d$ whose restriction to $E$ is linearly equivalent to $\o_E(\nu p)$ and let $n = d-g+1= \h^0(X,L)$ and $c=d-\nu=\deg_C(L)$. In our applications, $\nu=4$ and $L =\lbpow{\omega_X}{4}$, but for the moment using general $d$, $n$ and $\nu$ simplifies notation.

Consider the embedding of X in $\PP^{n-1}$ associated to $L$. It follows directly from Riemann-Roch that the linear spans $V_C$ of $C$ and $V_E$ of $E$ in $\PP^{n-1}$ are of dimensions $c-g+1 = n-\nu$ and $\nu-1$ respectively and that their intersection is $\{p\}$. Letting $l=n-\nu+1$, we can therefore choose homogeneous coordinates $x_1, \ldots, x_n$ such that $ x_1=  \ldots = x_{l-1} = 0$ defines $V_E$, $ x_{l+1}= \ldots = x_{n} = 0 $ defines $V_C$, and $p$ is the point where all the $x_i$ except $x_l$ vanish. 

For $ j \ge 1 $, we will confound $x_{l+j}$ with the section of $\HO{E}{\restrictedto{L}{E}}$ it determines and write $\ord_p(x_{l+j})$ for the order of vanishing at $p$ of this section. Again by Riemann-Roch, we may choose $x_{l+j}$ so that  $\ord_p(x_{l+j})=j$ for $1 \le j \le \nu-2$ and choose $x_n$ so that $\ord_p(x_n)=\nu$.

Define $\rho$ to be the $1$-ps subgroup of $\ssl(n)$ acting by $\diag(t^{r_1}, \cdots, t^{r_n})$ in these coordinates where
\begin{equation*}
r_i
=\begin{cases}
\nu & \text{if $i \le l$}\\
\nu-j & \text{if $i= l+j$ and $1 \le j \le \nu-2$}\\
0 & \text{if $i=n=l+\nu-1$}
\end{cases}
\end{equation*}
Note that, for $ j \ge 0$, this gives $x_{l+j}$ weight equal to $e-\ord_p(x_{l+j})$. 

Our sign convention is that the Hilbert-Mumford index $\mu([X]_m,\rho)$ is \emph{minus} the least weight with which $\rho$ acts on an $\thst{m}{th}$-Hilbert point of $(X,L)$ with respect to the canonical linearization of the  $\ssl(n)$-action. So our Numerical Criterion is
\begin{equation}\label{numericalcriterion}
\mu([X]_m,\rho) = -\bigl(w_{\rho}(m)- m P(m) \alpha(\rho)\bigr)
\end{equation}
where $w_{\rho}(m)$ is the least weight of any basis of $\HO{X}{\lbm{L}}$ consisting of the restrictions to $X$ of monomials of degree $m$ in the $x_i$, $\alpha(\rho)$ is the average \hbox{$\rho$-weight} of a coordinate on $\PP^{n-1}$ and $P(m)$ is the Hilbert polynomial of $(X,L)$.
Thus, with respect to $\rho$, $[X]_m$ is Hilbert stable, strictly semistable or unstable if and only if $\bigl(w_{\rho}(m)- m P(m) \alpha(\rho)\bigr)$ is, respectively, negative, $0$ or positive and hence $\mu([X]_m,\rho)$ has the opposite sign. Both terms are numerical (that is, integer \emph{valued}) polynomials of degree $2$ in $m$ and Chow stability for $(X,L)$ is similarly determined by the sign of the leading coefficient of the difference. For details, see, for example,~\cite{Moduli}*{Proposition 4.23 and Lemma 4.26} or \cite{MorrisonGIT}*{Proposition 3}).

In our examples, $P(m)=md-g+1$ by Riemann-Roch and an easy calculation shows that  
\begin{equation*}\alpha(\rho) = \frac{\nu n-\sum_{j=1}^{\nu-2} j - \nu}{n} = \nu- \frac{\nu^2-\nu+2}{2n}\,
\end{equation*}
 so that
\begin{equation}\label{mpmalpha}
m P(m) \alpha(\rho) = m^2\bigg[d\bigl( \nu- \frac{\nu^2-\nu+2}{2n}\bigr)\bigg]+	m\bigg[(1-g)\bigl( \nu- \frac{\nu^2-\nu+2}{2n}\bigr)\bigg]\,.
\end{equation}

\subsection*{Stability of elliptic tails}
The proof of Schubert's Lemma~3.1 applied to our setup, shows that the $m^2$ coefficient of $w_{\rho}(m)$ is at least $\bigl(d-\frac{\nu}{2}\bigr)\nu$ (although only $L=\lbpow{\omega_X}{3}$ is considered). The next lemma is the sharpening of this estimate to an exact evaluation of $w_{\rho}(m)$ needed to apply his argument to $4$-canonical models.
\begin{lemma}\label{lambdaweightlemma}
~~For $m \ge 2$, $\displaystyle{w_{\rho}(m) = m^2\bigg[\bigl(d-\frac{\nu}{2}\bigr)\nu\bigg] + m\bigg[\bigl(\frac{3}{2}-g\bigr)\nu\bigg] -1}$.
\end{lemma}
\begin{proof}
For concision, we will henceforth understand all monomials to have degree $m$ and view them directly as sections of $\lbm{L}$ over $X$ or $E$ (eliding ``the restriction to''). Weights will always be $\rho$-weights. 

Let $W_r$ be the span in $\HO{X}{\lbm{L}}$ of all monomials of weight at most $r$ and let $s=m\nu-r$.   We claim that
\begin{equation*}
	\dim(W_r)
	=\begin{cases}
	md-g+1 & \text{if $r =  m\nu$}\\
	r & \text{if $2 \le r \le  m\nu-1$}\\
	1 & \text{if $r=0$ or $r=1$}
	\end{cases}		
\end{equation*}
Given this, the lemma follows by elementary manipulations since the weight of any monomial basis of minimal weight is simply the sum of $r \big(\dim(W_r)-\dim(W_{r-1})\bigr)$ over $r$.

The first case in the claim is immediate from Riemann-Roch for $\lbm{L}$ on $X$. The others follow from the equality
\begin{equation}\label{weightspaces}
	W_{r} = \HO{E}{\restrictedto{\lbm{L}}{E}(-sp)} \text{~~for $r=0$ \text{and for} $2\le r \le m\nu-1$\,.}
\end{equation}
Since $\restrictedto{\lbm{L}}{E}(-sp)\cong \o_E(rp)$---recall that $\restrictedto{L}{E} \cong \o_E(\nu p)$, Riemann-Roch on $E$ implies that $\hO{E}{\restrictedto{\lbm{L}}{E}(-sp)}=r$.

If any monomial has weight $r=m\nu-s$ then it contains one or more factors $x_{l+j}$ with $j>0$ and hence vanishes on $C$. By construction, $s$ equals the sum of the orders of vanishing at $p$ of the factors of this type, hence $W_{r} \subset \HO{E}{\restrictedto{\lbm{L}}{E}(-sp)}$. 

If we next set $M_0 = x_n^m$, then $B_0 := \{M_0\}$ is a basis of $\HO{E}{\restrictedto{\lbm{L}}{E}(-m\nu p)}$ lying in $W_0$. Finally, for $r=2, \ldots, m\nu-1$, let $M_r$  be any monomial 
$$M_r := \prod_{k=1}^m x_{l+j_k} \text{~s.t. each $j_k \ge 0$ and} \sum_{k=1}^m j_k = s=m\nu-r\,.$$
Then, $M_r$ vanishes on $C$ because  some $j_k>0$.  By construction, $M_r$ has weight exactly $r$ and, since $x_l$ is non-zero at $p$, $M_r$ vanishes to order exactly $s$ at $p$. Thus, $B_r:= \{M_0, M_2, M_3, \ldots, M_r\}$ is a subset of $W_r \cap \HO{E}{\restrictedto{\lbm{L}}{E}(-sp)}$ of cardinality equal to $\hO{E}{\restrictedto{\lbm{L}}{E}(-sp)}$. But all the elements of $B_r$ except $M_r$ lie in $W_{r-1}$. By induction, $B_r$ is linearly independent and hence is a basis of $\HO{E}{\restrictedto{\lbm{L}}{E}(-sp)}$ which therefore lies in~$W_{r}$.
\end{proof}

We now specialize to the case where $\nu=4$ and $L=\lbpow{\omega_X}{4}$ so that $d= 8(g-1)$ and $n=7(g-1)$. From Lemma~\ref{lambdaweightlemma} we obtain  $w_{\rho}(m)= (32g-40)m^2 + (-4g+6)m - 1$ and, from~\eqref{mpmalpha}, $m P(m) \alpha(\rho) = (32g-40)m^2 + (-4g+5)m$.

\begin{corollary}\label{fourcanonicaltailsunstable}
If $X$ is a $4$-canonically embedded Deligne-Mumford stable curve with an elliptic tail then $\mu([X]_m, \rho) = -m+1$. Hence, with respect to $\rho$, such $X$ are Chow strictly semistable and $m$-Hilbert unstable for any $m \ge 2$.
\end{corollary}

There are many other examples with the last properties. For any $\nu$ we can take $\alpha = \frac{\nu^2}{\nu-2}$, $d=\alpha(g-1)$ and $n=(\alpha-1)(g-1)$ for $g-1$ divisible by $\nu-2$. When $\nu=3$, $g$ can be arbitrary. We have checked that none of these Chow strictly semistable examples are Hilbert semistable. 

Lemma~\ref{lambdaweightlemma} also reproves Schubert's Lemma~3.1 by taking $\nu=3$ and $L=\lbpow{\omega_X}{3}$.

\subsection*{Stability of cuspidal tails}

In this section, we consider only $4$-canonical curves of genus $g$ (so $d=8(g-1)$, $n=7(g-1)$ and $l=7g-10$). We want to consider a curve $Y$ with a cuspidal tail: that is, $Y$ is the join of a curve $C$ of genus $g-1$ and a rational cuspidal curve $R$ (of arithmetic genus $1$) meeting $C$ in a single node at~$p$. 

The key remark is that the definition of $\rho$ extends unchanged to these curves. Once again $\restrictedto{\omega_Y}{R} = \o_R(p)$, we may choose coordinates such that  $ x_1=  \ldots = x_{l-1} = 0$ defines $V_R$, $ x_{l+1}= \ldots = x_{n} = 0 $ defines $V_C$, and $p$ is the point where all the $x_i$ except $x_l$ vanish. Moreover, if we parameterize $R$ by \[
[s,t] \mapsto [0, \ldots, 0, t^4, st^3, s^2t^2, s^4].
\]
then $\ord_p(x_{l+j}) = j$ for $j\ge 1$ with its $\rho$-weight is the power $n-l+j$ of $t$. Finally, $\rho$ stabilizes $R$, and, as in Corollary~\ref{fourcanonicaltailsunstable}:
  
\begin{lemma}\label{L:Cstar} If $Y$ is a $4$-canonical curve with a cuspidal tail, then $\mu([Y]_m, \rho) = m - 1$. Hence, with respect to $\rho$, such $Y$ are Chow strictly semistable and $m$-Hilbert unstable for all $m\ge 2$. 
\end{lemma}
\begin{proof} The approach of Lemma~\ref{lambdaweightlemma} works here but there's a more elementary and concrete alternative. Because $4$-canonical curves are embedded by complete linear series with no higher cohomology as $2$-regular subschemes in the sense of Castelnuovo-Mumford, we can apply Proposition~3.17 of~\cite{HassettHyeonFlip}. It asserts that $w_F(m) = (m-1) (am+b)$ for all $m\ge 2$ and hence is determined by $w_F(2)$ and $w_F(3)$ (by an explicit formula that we will not need). 
	 
The degree two monomials not in the initial ideal $\init_\rho(R)$ of $R$ must pullback to $s^{8-i}t^{i}$ (of weight $i$) where $i$ runs from $0$ to $8$ skipping $1$. Hence their weights sum to $35$.
On the other hand,
$\rho$ acts on $D$ with constant weight $4$ and
the number of monomials not in the initial ideal of $D$ vanishing at the
point of attachment $p$ is $h^0(L^{\otimes 2}(-p)) = 15g-23$
where $L = \omega_C^{\otimes 4}(4p)$. Hence the total weight in degree $2$ is
$
2\cdot 4\cdot (15g - 23) + 35 = 120 g - 149
$.
Likewise, the degree three monomials not in $\init_\rho(R)$
pull back to $s^{12-i}t^{i}$ with $0\le i\le 12$ and $i \neq 1$
and have weights summing to~$77$. Adding the weight $3\cdot 4\cdot h^0(D, L^{\otimes 3}(-p)) = 12 (23g-35)$ of any monomial basis of sections on $D$
vanishing at $p$, we get a total weight in degree $3$ of $276g-343$. 

On the other hand,~\eqref{mpmalpha} gives  $mP(m)\alpha_{\rho}= m(8m-1)(4g-5)$ which yields $120g-150 $ for $m=2$ and $276g-345 $ for $m=3$. This shows that $\mu([Y]_m, \rho) = -m+1$ for $m=2$ and $3$ and hence for all $m$.
\end{proof} 

Let $G$ be a reductive group, $\h \subset \mathbb \PP(V)$ be a projective
variety with a linear $G$ action. Let    $\rho:\mathbb G_m \to G $ be a
one parameter subgroup with fixed points $\h^{\rho}$.
For each $x^{\star} \in \h^{\rho}$, the {\em $\rho$-basin of attraction}
is defined
\[
A_\rho(x^\star) := \left\{ x\in \h \, | \, \lim_{t \to 0} \rho(t)\cdot x = x^\star
\right\}.
\]

If $x \in A_\rho(x^\star)$ and $\mu(x^\star, \rho) = 0$,
then $x$ is semistable if and only if $x^\star$ is.  Since the orbit closure
of any $x \in A_\rho(x^\star)$ contains $x^\star$, if $\mu(x^\star, \rho) = 0$ and $x^\star$ is semistable,
then $x$ and $x^\star$ are identified in the quotient space $\h/\!\!/G$. We want to apply this to the Chow point $[Y]$ of the curve with a cuspidal tail.

\capdraw{F:1n1t}{Basin of attraction of $Y$}{.5}{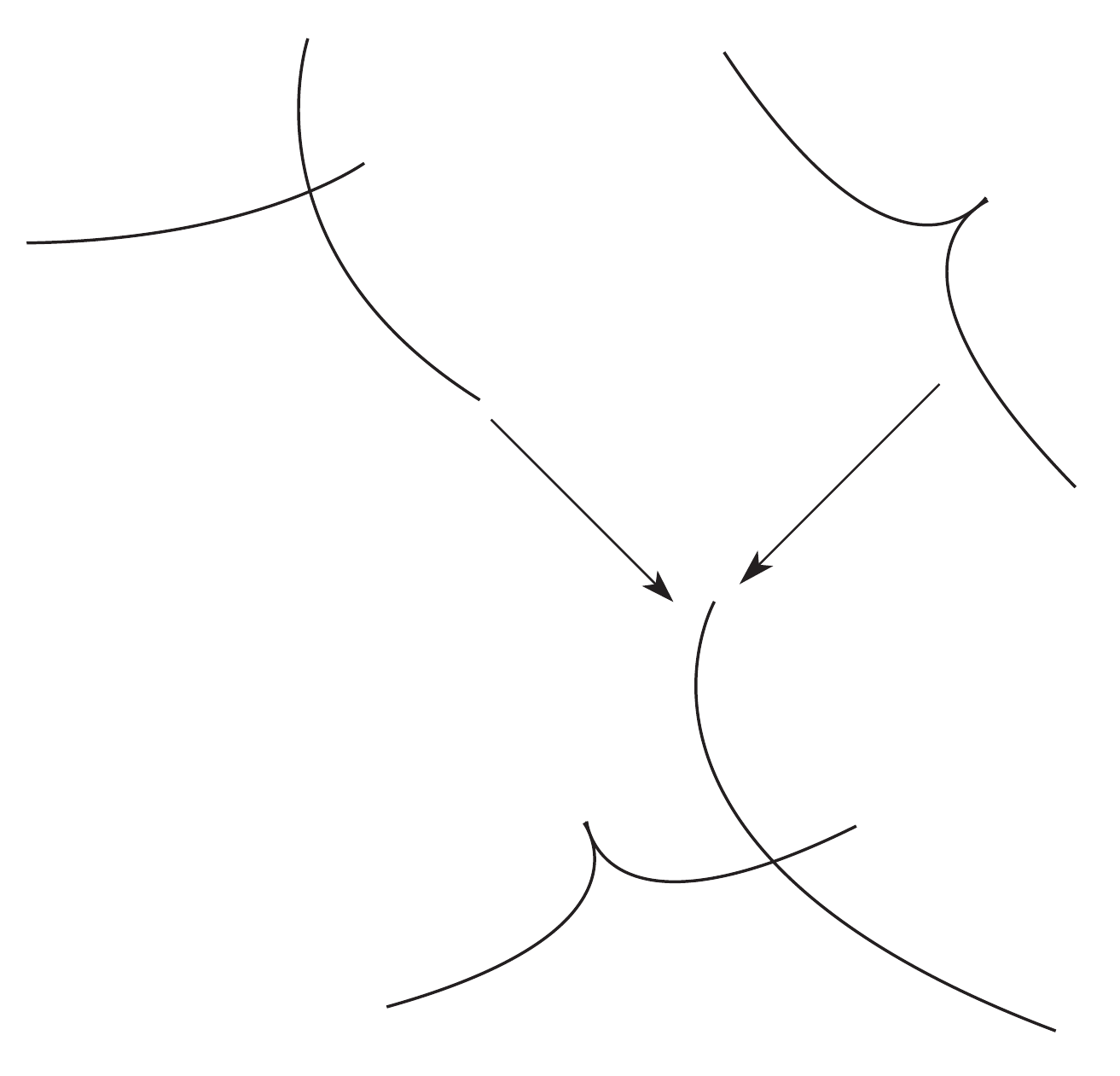}{%
\etext{.17}{.90}{$X$}
\etext{.55}{.35}{$Y$}
\wtext{.80}{.90}{$Z$}
\etext{.33}{.07}{$R$}
\stext{.05}{.79}{$E$}
\nwtext{.30}{.82}{$p$}
\ntext{.69}{.18}{$p$}
\netext{.49}{.54}{$\rho$}
\nwtext{.75}{.54}{$\rho^{-1}$}
\ntext{.12}{.75}{elliptic tail}
\wtext{.38}{.70}{$C$}
\wtext{.82}{.13}{$C$}
\etext{.47}{.16}{rational cuspidal tail}
\stext{.52}{.26}{$q$}
\wtext{.95}{.74}{$C$ pinched to}
\wtext{.96}{.68}{an ordinary}
\wtext{.98}{.62}{cusp at $q$}
\swtext{.90	}{.82}{$q$}
}%

\begin{lemma}\label{L:basin} The basin of attraction $A_\rho([Y])$ 
\resp{$A_{\rho^{-1}}([Y])$} contains
smoothings of the cusp \resp{node} but no smoothings of the node \resp{cusp}.
\end{lemma}
\begin{proof} 
At the cusp $q=[0,\dots,0,1]$, the deformation
space is given by $y^2 = x^3 + a x + b$, and $\rho$ acts on $x$ (which 
may be taken to be (the image of) $\frac{x_{n-1}}{x_n}$ in the completed local ring) with weight $2$. Hence $\rho$ acts on $a$ and $b$ with weights $4$ and $6$, and the basin of attraction
contains smoothings of the cusp. On the other hand, local coordinates at the node are given by dividing by $x_{n-3}$. In these, $\rho$ acts on 
the tangent line to $E$ at the node $p$ (spanned by $\frac{x_{n-2}}{x_{n-3}}$) with weight $-1$ and on the tangent line to $C$ at $p$ trivially. Hence the basin of attraction does not admit a smoothing of the node. Since $\rho^{-1}$ acts on the deformation spaces with  
weights of the opposite signs, the opposite conclusions hold. 
\end{proof} 

This reproves Corollary~\ref{fourcanonicaltailsunstable}. As a complement, we sketch a direct proof that 
$Z$ is $m$-Hilbert stable and Chow strictly semistable with respect to $\rho^{-1}$ which we normalize to have weights $[0,\ldots, 0, 1, 2, 4]$ with $x_1=\ldots=x_k=0$ for $k=n-1, n-2$ and $n-3$ defining the cusp $q$, its tangent line and its osculating plane respectively. Thus $x_n, x_{n-1}, x_{n-2}$ and $ x_{n-3}$ vanish, respectively, to orders $0,1,2$ and $4$ at $q$.

\begin{lemma}\label{cusplambdainverse}
If $Z$ is a $4$-canonical curve with a cusp, then $\mu([Z]_m, \rho^{-1}) = m-1$. Hence, with respect to $\rho^{-1}$, such $Z$ are Chow strictly semistable and $m$-Hilbert stable for all $m\ge 2$. 
\end{lemma}
\begin{proof}
For $j \ge 0$, $x_{l+j}$ has $\rho^{-1}$-weight equal to $4-\ord_q(x_{l+j})$. This implies, in the notation of Lemma~\ref{lambdaweightlemma} (for $\rho^{-1}$), that $W_{4m-s}= \HO{Z}{\lbpow{\omega_Z}{4m}(-sq)}$ for $0< s \le 4m$. Thus, $\dim(W_r/W_{r-1})=1$ for $0<r<4m-2$ and $\dim(W_{4m}/W_{4m-2})=1$ and the least weight of a basis is $\bigl(\sum_{r=1}^{4m} r\bigr) - (4m-1) = 8m^2-2m+1$. Since the weights total~$7$, $mP(m)\alpha_{\rho^{-1}}= m(8m-1)(g-1)\frac{7}{7(g-1)}=8m^2-m$, which gives the claimed~$\mu$. 
\end{proof}

\begin{proposition}\label{P:chowfour}
The $4$-canonical Chow points of 
$X$, $Y$ and $Z$ are Chow strictly semistable.
\end{proposition}
\begin{proof}
In view of Corollary~\ref{fourcanonicaltailsunstable} and Lemmas~\ref{L:Cstar} and~\ref{cusplambdainverse}, we must eliminate the possibility that $Y$ is Chow unstable. By Lemma~\ref{L:basin}, we can do this by showing that any of the three curves is Chow semistable. 

Let $B = \spec(k[[t]])$ with generic point $\eta$ and closed point $0$ and let $S \to B$ be a family of curves over $B$ with generic fiber smooth of genus $g$ and special fiber $S_0$ isomorphic to $Y$. 
After making a base change, if necessary, we can assume that there is a family $S'\to B$ with general fiber $S'_{\eta}$ isomorphic to $S_{\eta}$ and with Chow semistable special fiber $S'_0$. 

Schubert's arguments~\cite{Schubert}*{Section 2} show that $S'_0$ can have only nodes and cusps as singularities and that it cannot contain any semistable chains of smooth rational curves. 
By stable reduction, on the other hand, $S'_0$ must have the same Deligne-Mumford stabilization as $Y$. These two conditions force $S'_0$ to be one of the curves $X$, $Y$ or $Z$.
\end{proof}

\subsection*{The $4$-canonical GIT quotient spaces} Let $\Hilb_{g,\nu}$ denote the closure of the locus of smooth $\nu$-canonically embedded curves in the Hilbert scheme of curves of genus $g$ and degree $2\nu(g-1)$ in $\PP^{(2\nu-1)(g-1)-1}$ and let $\Chow_{g,\nu}$ denote the corresponding locus in the Chow variety.  A complete connected curve is {\it weakly pseudostable} if it has only nodes and ordinary cusps as singularities and any smooth rational component meets the rest of the  curve in at least three points. 

Combining Schubert's arguments with Corollary~\ref{fourcanonicaltailsunstable}, Lemmas~\ref{L:Cstar} and~\ref{cusplambdainverse}, and Proposition~\ref{P:chowfour} gives: 
\begin{enumerate}
\item In ${\Hilb_{g,4}}$, the stable locus parameterizes pseudostable curves, and there are no strictly semistable points. 
\item In ${\Chow_{g,4}}$ the semistable locus parameterizes weakly pseudostable curves and the stable locus Deligne-Mumford stable curves without elliptic tails. 
\item Two weakly pseudostable curves are identified in $\chowq$ if and only if they have a cusp or an elliptic tail, and their Deligne-Mumford stabilizations are the same (see Figure~\ref{F:1n1t}).
\end{enumerate}

As noted in the introduction, the arguments in~\cite{Schubert} go through without significant change to show that $\hilbq \simeq  \mpsg$. The fundamental cycle map $FC$ from the Hilbert scheme to the Chow variety descends to a map 
$$\varpi: \hilbq \to \chowq \,.$$
By our classification of semistable points, $\varpi$ is an isomorphism away from the locus of elliptic tails, is everywhere injective and induces a bijection on closed points.  

The proof of~\cite{HassettHyeonFlip}*{Corollary~6.3} shows that, for any $\nu \ge 2$, the $\nu$-canonical locus $\Hilb_{g,\nu}$ is smooth at points parameterizing weakly pseudostable curves. In particular, $\semi{\Hilb_{g,4}}$  is smooth and its good quotient $\hilbq$ is normal. 

Koll{\'a}r~\cite{KollarRational}*{Section~I.6} implies that there is a flat universal family of cycles over $\Chow_{g,\nu}$ that yields a regular inverse morphism to $FC$ over this locus and his Theorem~7.3.1 then shows that $\chowq$ is also normal. Since it is an injective birational projective morphism between normal varieties, $\varpi$~is an isomorphism.
We have shown that:

\begin{theorem}\label{T:main} $
	\chowq \simeq \hilbq \simeq \mpsg$.
\end{theorem}

\medskip

\ifbells

\fi

\end{document}